\documentclass[12pt]{article}
\usepackage{float}

\usepackage{setspace}
\usepackage{amsmath}
\usepackage{amssymb, amsmath, amsthm}
\usepackage{mathrsfs}
\usepackage{tabulary}
\usepackage{rotating}
\usepackage[para]{threeparttable} 
\usepackage{enumerate}
\usepackage{setspace}
\usepackage{multirow}
\usepackage{paralist}
\usepackage{listings}
\usepackage[T1]{fontenc}
\usepackage[utf8]{inputenc}
\usepackage[english]{babel}
\usepackage{parskip}
\usepackage{dsfont}
\usepackage{hyperref}
\usepackage{subcaption}
\usepackage{algorithm}
\usepackage{algpseudocode}
\usepackage{dsfont}
\usepackage{arydshln}
\usepackage{multirow}

\usepackage{graphicx}
\usepackage{url}
\usepackage{booktabs}
\usepackage{enumitem}
\newlist{subquestion}{enumerate}{1}
\setlist[subquestion,1]{label=(\alph*)}
\usepackage{pdflscape}
\usepackage{float}
\usepackage{adjustbox}
\usepackage{etaremune} 

\usepackage{authblk}
\usepackage{etoolbox}
\usepackage{longtable}
\usepackage[numbers,super,sort&compress]{natbib}
\usepackage[usenames,dvipsnames]{color}
\usepackage{natbib}
\usepackage{pdfpages}

\usepackage{caption}
\usepackage{hyperref}

\usepackage[hmargin=2.54cm,vmargin=2.54cm]{geometry}
\usepackage[utf8]{inputenc}
\usepackage[hang,flushmargin]{footmisc}

\DeclareMathOperator*{\argmin}{arg\,min}

\newcommand{\E}{\mathds{E}}

\newcommand{\cd}{\,}

\newcommand{\m}{\mathsf{m}}
 
\newcommand{\Q}{\mathsf{Q}}

\newcommand{\one}{\mathds{1}}
\newcommand{\rr}{\alpha(X)}

\newtheoremstyle{customnote}%
  {\topsep}
  {\topsep}
  {\upshape}
  {}
  {\bfseries}
  {.}
  {.5em}
  {\thmname{#1}\thmnumber{ #2}\thmnote{ \normalfont\itshape #3}}

\theoremstyle{customnote}
\newtheorem{theorem}{Theorem}
\newtheorem{example}{Example}

\theoremstyle{definition}
\newtheorem{definition}{Definition}

\newtheorem{alg}{Algorithm}

\newenvironment{keywords}
{\bgroup\leftskip 20pt\rightskip 20pt \small\noindent{\bf Keywords:} }%
{\par\egroup\vskip 0.25ex}

\begin{document}

\title{Riesz representers for the rest of us}

\vspace{3cm}

\author[1]{Nicholas T. Williams}
\author[1]{Oliver J. Hines}
\author[1]{Kara E. Rudolph}
\affil[1]{\footnotesize Department of Epidemiology, Mailman School of Public Health, Columbia University, New York, NY}
\date{}

\maketitle

\begin{abstract}
    The application of semiparametric efficient estimators, particularly those that leverage machine learning, is rapidly expanding within epidemiology and causal inference. This literature is increasingly invoking the Riesz representation theorem and Riesz regression. This paper aims to introduce the Riesz representation theorem to an epidemiologic audience, explaining what it is and why it's useful, using step-by-step worked examples.
\end{abstract}


\section{Introduction}\label{section:intro}
The use of terms like ``Riesz regression/learning" and ``Riesz representers" in the causal inference literature has grown rapidly. Whereas only a handful of scientific articles invoked these terms 10 years ago, nearly 50 such articles used these terms in the past year alone. Epidemiologists may wonder why so many of their colleagues in (bio)statistics are suddenly including the above phrases in their papers and talks, and may ask themselves: What is a Riesz representer? What is Riesz regression? Why should I care? 
With these questions in mind, we offer a gentle introduction for epidemiologists.

First, we provide some background on \textit{doubly robust} (DR)-estimators. Over the past 30 years, epidemiological analyses have increasingly used DR-estimators that are rooted in semiparametric statistical theory for the estimation of causal effects.\citep{smith2023application} Initial theoretical work focused on estimating the Average Treatment Effect (ATE), with the proposal of the augmented inverse probability weighted (AIPW) estimator,\citep{robins_estimation_1994} which can be viewed as a combination of g-computation \citep{robins1986new} and inverse probability weighting (IPW).\citep{HorvitzThompson1952} It was later noted that the AIPW estimator is asymptotically consistent---meaning, as sample size grows, the estimator converges to the true value---provided that either the outcome model or the propensity score model is correctly specified,\citep{scharfstein_adjusting_1999} a property that would become known as double robustness.\citep{bang_doubly_2005,van2006targeted,funk_doubly_2011} DR-estimators have an additional property, known as \textit{rate double robustness}, which allows for leveraging machine learning for estimation of the outcome and exposure models whilst still being able to construct confidence intervals and perform statistical hypothesis testing. For this reason, DR-estimators are sometimes referred to as \textit{debiased} since they remove biases which arise from the naive use of machine learning in estimation. Similar estimators have since been developed for a wide range of scenarios; select examples include: mediation analysis,\citep{tchetgen2012semiparametric,zheng2012targeted,diaz2021nonparametric} continuous exposures, \citep{munoz2012population,haneuse2013estimation,kennedy2017non} treatment effect heterogeneity,\citep{levy_fundamental_2021} survival outcomes,\citep{hubbard2000nonparametric,benkeser2018improved,cai2020one} and longitudinal data.\citep{van2012targeted,diaz_nonparametric_2023} 

A central task in constructing DR-estimators is derivation of the efficient influence function (EIF) of the \textit{statistical estimand}, where we refer to statistical estimand as the quantity to be estimated that is a function of observed data (distinguishing it from a \textit{causal estimand}, which is a quantity that is a function of counterfactual variables). Many statistical estimands in the causal inference literature can be expressed as a bounded linear functionals composed of iterated conditional expectations (we discuss this in detail in \S\ref{sec:linearfunc}). A linear functional is a linear mapping that takes as input one or more functions and returns a real number. In such cases, the EIF is composed of a conditional outcome regression and a quantity called the \textit{Riesz representer}---both of which require estimation. Many epidemiologists will be familiar with the Riesz representers of several canonical estimands in epidemiology. For example, we will show that the Riesz representers of the ATE and the average treatment effect on the treated (ATT) recover the IPW expressions used to construct Horvitz-Thompson style reweighting estimators.\citep{HorvitzThompson1952} In fact, the existence of these IPW weights is a direct consequence of the Riesz representation theorem.

The standard approach for estimating IPW weights is to fit a propensity score model and invert the subsequent propensity score predictions. However, this may lead to unstable behavior in finite samples. For example, one may observe extreme IPW weights when practical overlap violations between the treated and untreated populations occur.\citep{kang_demystifying_2007,damour_overlap_2021} As such, there has been substantial interest in estimating Riesz representers directly from data, without needing to estimate constituent components (such as propensity scores). The main proposals to achieve this are to use \textit{balancing weights}, which are weights with minimal dispersion that approximately satisfy known properties of the Riesz representer\citep{hainmueller_entropy_2012,imai_covariate_2014,zubizarreta_stable_2015,ben-michael_balancing_2021}, and \textit{Riesz regression}, which uses empirical risk minimization to learn Riesz representers using a bespoke loss function.\citep{chernozhukov2021automatic,chernozhukov2022automatic,chernozhukov2022riesznet,susmann2024longitudinal,hines2024automatic,van2024automatic,liu2024general,lee2025rieszboost,hines_learning_2025} These approaches have the added benefit that they allow for automatic debiasing, in the sense that one can construct DR-estimators without knowing the full functional form of the estimand's EIF. This may be useful if the EIF is difficult to derive using conventional methods.\citep{kennedy2024semiparametric,Hines_Dukes_Diaz-Ordaz_Vansteelandt_2022,renson2025pulling} (We note, however, that the benefit of automatic debiasing is arguably more theoretical than practical, since we are unaware of any DR-estimators for which the EIF has not already been derived.)

In this paper, through worked examples, we aim to provide an accessible introduction to the Riesz representation theorem and explain how Riesz representers relate to DR-estimation. Additionally, we introduce a simple algorithm for deriving the functional form of the EIF based on Riesz representers. Lastly, we review how one can estimate Riesz representers via Riesz regression. 

\section{Preliminaries}

This tutorial is concerned with estimation, which is one step in the process of answering an epidemiologic question. First, epidemiologists use their subject-matter knowledge to propose a causal model (i.e., Neyman-Rubin model,\citep{rubin1974estimating} non-parametric structural equation model,\citep{Pearl2000Causality} directed acyclic graph, or single-world intervention graph\citep{richardson2013single}) and a causal estimand, which is the causal quantity of interest (e.g., ATE), defined using unobserved counterfactual variables. Then, we identify the assumptions needed to equate this causal estimand with a statistical estimand that is a function of only observed data. This step-by-step process of showing the assumptions needed to equate the causal estimand with the statistical estimand is called identification. After identification, we are ready to estimate this statistical estimand. Throughout, we will assume that each example causal estimand can be identified from the observed data by a statistical estimand.

Let $X = (W, A, M, Y)$ denote the observed data. Here, $A$ represents a treatment or exposure of interest, $Y$ an outcome, $M$ a vector of mediators between $A$ and $Y$, and $W$ a vector of confounders of the effect of $A$ on $Y$. Random variables will be denoted using upper-case and realizations of a random variable will be denoted using lower-case (e.g., $x$ is a realization of $X$). Throughout, we use $f(\cdot)$ to denote the density or probability mass function of a random variable (e.g., $f(y\mid a,w)$ is the density of $Y=y$ conditional on $A=a$ and $W=w$ and $f(A=1)$ is the probability that $A=1$). We remind the reader that a conditional expectation of a continuous random variable is defined using an integral (e.g., $\E[Y\mid A=a,W=w] = \int y\,f(y\mid a, w)\,dy$). 

\section{What is the Riesz representation theorem?}

\subsection{Estimands as bounded linear functionals}\label{sec:linearfunc} 

We focus on estimands that can be defined as \textit{bounded linear functionals} of a sequence of linear operations and conditional outcome expectations. As noted above, a functional is a function that takes as input other functions and returns a real number. For example, assuming $A$ is binary, $\E\{\E[Y\mid A=a,W]\}$ is a bounded linear functional of $\E[Y\mid A=a, W=w]$. First, we regress $Y$ on $A$ and $W$. Then, we evaluate predictions of that regression at $A=a$ and random variable $W$, which is a linear mapping. Finally, we take the marginal expectation of that output. Intuitively, bounded linear functionals are important because if an estimator for the conditional expectation converges to its true value with increasing sample size, then the output of the functional applied to the estimator will also converge to the true value. We provide a mathematical definition of boundedness and linearity in the eAppendix.

 In what follows, we denote outcome regressions as $\Q_k(X)$, and linear mappings of regressions as $\m_k(\Q_k;x)$ to highlight that $\m_k$ is a mapping of the regression $\Q_k$ dependent on a realization $x$. Here, $k$ indexes the regression (this will become clear in the examples that follow). For example, if $\Q_k(a,w) = \E[Y\mid A=a,W=w]$, the mapping $\m_k(\Q_k;w) = \Q_k(1,w)$ evaluates the regression setting $A=1$ at different values of $w$. Note that $\m_1$ may simply be the identity functional that returns its input unchanged (i.e., $\m_1(\Q_1;\cdot) = \Q_1(\cdot)$) and will always represent the outermost operation in the composition defining the target functional. This setup allows us to characterize many of the common causal effects in the literature, and will allow us to introduce the Riesz representation theorem.

\begin{example}[Mean outcome among the treated]\label{rct}
    In a two-arm randomized controlled trial (RCT), the mean outcome in the trial population had all units been assigned treatment is identified from the observed data as the statistical estimand $\theta = \E[Y\mid A=1]$. The conditional outcome regression is denoted $\Q_1(a) = \E[Y\mid A=a]$. In this very simple example, the functional is denoted $\m_1(\Q_1) = \Q_1(1)$, meaning the prediction from the outcome regression setting $A=1$. Using the previously described notation, $\theta = \m_1(\Q_1)$.
\end{example}

\begin{example}[Average treatment effect]\label{ate}
    The ATE is defined as the difference in counterfactual outcomes had all units in the population received treatment versus if no units had received treatment. Under suitable conditions\citep{Hernan2020_Observational}, the ATE is identified from the observed data as
    \begin{equation}\label{eq:ate}
        \theta = \E\{\E[Y\mid A=1,W] - \E[Y\mid A=0,W]\}. 
    \end{equation}
    Let $\Q_2(a, w) = \E[Y\mid A=a,W=w]$ and $\m_2(\Q_2;W) = \Q_2(1,W) - \Q_2(0,W)$. Additionally, let $\Q_1 = \E[\m_2(\Q_2;W)]$ and $\m_1(\Q_1) = \Q_1$ (identity). 
Again, $\theta = \m_1(\Q_1)$.
\end{example}

\begin{example}[Mean outcome under control for the treated]\label{att}
    The counterfactual mean outcome under control for the treated is the expected value of the outcome had all units who received treatment instead received the control; importantly, it's one component of the ATT. Under suitable conditions\citep{Hernan2020_Modification} similar to that of the ATE, the mean outcome under control for the treated is identified from the observed data as
    \begin{equation}\label{eq:att}
        \theta = \E\{\E[Y\mid A=0,W]\mid A= 1\}.
    \end{equation}
    Let $\Q_2(a,w) = \E[Y\mid A=a,W=w]$ and $\m_2(\Q_2;W) = \Q_2(0,W)$. Additionally, let $\Q_1(a) = \E[\m_2(\Q_2;W)\mid A = a]$ and $\m_1(\Q_1) = \Q_1(1)$. $\theta = \m_1(\Q_1)$.
\end{example}

\begin{example}[Natural direct effect]\label{nde}
    Mediation analysis is concerned with decomposing the effect of an exposure on an outcome via post-exposure variables $M$, referred to as mediators. For example, the direct effect refers to the effect of an exposure via pathways that do not involve the mediator and can be measured through a variety of estimands. One such estimand is the natural direct effect (NDE)\citep{pearl2022direct}, which is the expected difference in the counterfactual outcome had all units been exposed versus unexposed, and the mediator takes the value it would naturally take under no exposure. Again, under suitable conditions\citep{pearl2022direct}, the NDE of a binary exposure is identified from the observed data as a contrast of the statistical estimand 
    \begin{equation}\label{eq:nde}
        \theta(a') = \E[\E\{\E[Y\mid A=a', M, W]\mid A=0, W\}]
    \end{equation}
    setting $a'=1$ versus $a'=0$.
    Let $\Q_3(a,m,w) = \E[Y\mid A=a, M=m, W=w]$ and $\m_3(\Q_3;M,W) = \Q_3(a',M,W)$. Additionally, let $\Q_2(a,w) = \E[\m_3(\Q_3;M,W)\mid A=a,W=w]$ and $\m_2(\Q_2;W) = \Q_2(0, W)$, and $\Q_1 = \E[\m_2(\Q_2;W)]$. Then $\theta(a') = \m_1(\Q_1) = \Q_1$.
\end{example}

\subsection{An alternate expression using Riesz representers}

In \S\ref{sec:linearfunc}, we expressed each of the example statistical estimands as iterative compositions of linear mappings and conditional expectation. We now show that each estimand can also be expressed as an inner product of a unique set of weights and the outcome $Y$. These two general approaches of expressing the statistical estimand---iterative conditional expectations and weights multiplied by the outcome---are well-known and correspond to identification using a g-computation formula vs. identification using an IPW formula.\citep{HernanRobins2020} In fact, \textit{using the IPW formula for identification is a direct consequence of the Riesz representation theorem.} For example, using the IPW formula for identification, the mean outcome among the treated (Example \ref{rct}) can alternatively be written as
\begin{equation}\label{eq:ipw:rct}
    \begin{split}
        \theta &= \E[Y\mid A=1] \\
            &= \int y \cd f(y\mid A=1)\,dy \\
            &= \frac{1}{f(A=1)}\int y \cd f(y, A=1)\,dy\\
            &= \frac{1}{f(A=1)} \sum_{a=0}^1 \int y \one(a=1) f(y, a) \,dy \\
            &= \E\bigg[\frac{\one(A=1)}{f(A=1)}\cd Y\bigg] \\
            &= \E\bigg[\frac{\one(A=1)}{f(A=1)}\cd \E[Y\mid A] \bigg].
    \end{split}
\end{equation}
The second equality follows from the definition of a conditional expectation, the third from the definition of conditional probability and because $f(A=1)^{-1}$ doesn't depend on $y$. The fourth equality is obtained by noticing that $\int y \cd f(y, A=1)\,dy = \sum_{a=0}^1 \int \one(a = 1)\cd y\cd f(y, a)\,dy$, where $\one(a=1)$ is the indicator function that returns $1$ when $a=1$ and $0$ when $a=0$. The last equality follows from the law of total expectation. This re-expression is a direct consequence of the Riesz representation theorem. 

Similarly, the IPW identification result of the ATE\citep{rosenbaum1983central} (Example \ref{ate}) is given by
    \begin{equation}\label{eq:ipw:ate}
        \begin{split}
            \theta &= \E\{\E[Y\mid A=1,W] - \E[Y\mid A=0,W]\} \\
                &= \iint y\cd\big\{f(y\mid A=1,w) - f(y\mid A=0,w)\big\} f(w)\,dy\,dw\\
                &= \iint y \cd \bigg\{\frac{f(y,A=1,w)}{f(A=1,w)} - \frac{f(y,A=0,w)}{f(A=0,w)}\bigg\}f(w)\,dy\,dw\\
                &= \iint y \cd \bigg\{\frac{f(y,A=1,w)}{f(A=1\mid w)} - \frac{f(y,A=0,w)}{f(A=0\mid w)}\bigg\}\,dy\,dw\\
                &= \E\bigg[\bigg\{\frac{\one(A=1)}{f(A=1\mid W)} - \frac{\one(A=0)}{f(A=0\mid W)}\bigg\}\cd \E[Y\mid A,W] \bigg],
        \end{split}
    \end{equation}
where each equality follows from analogous operations shown in Eq. \ref{eq:ipw:rct}.

\begin{theorem}[Riesz representation theorem---Informal]\label{riesz}
    Every bounded linear functional $\theta = \m_1(\Q_1)$ (where $\Q_1$ is the final iterated expectation) can be expressed as an inner product of the initial outcome regression $\Q_K$ (or the outcome of that regression e.g., $Y$) and weights $\rr$.
    \begin{equation}
       \theta = \m_1(\Q_1) =\E[\rr\cd Y] = \E[\rr\cd\Q_K(X)]  
    \end{equation}
\end{theorem}
The weights, $\alpha(X)$, in Theorem \ref{riesz} are what are referred to as Riesz representers in the double and de-biased machine learning literature.\citep{chernozhukov2022automatic}

The Riesz representer for Example \ref{rct} is $\alpha(A) = \one(A=1) f(A=1)^{-1}$; the Riesz representer for Example \ref{ate} is $\alpha(A,W) = \one(A=1)f(A=1\mid W)^{-1} - \one(A=0)f(A=0\mid W)^{-1}$. Let's calculate the Riesz representers for the remaining example estimands. For the mean outcome under control for the treated (Example \ref{att}), we have
\begin{equation}\label{eq:ipw:att}
    \begin{split}
      \theta &= \E\{\E[Y\mid A=0,W]\mid A=1\}\\
        &= \iint y\cd f(y\mid A=0,w)\cd f(w\mid A=1)\,dy\,dw\\
        &= \frac{1}{f(A=1)}\iint y\cd f(y\mid A=0,w)f(A=1,w)\,dy\,dw\\
        &= \frac{1}{f(A=1)}\iint \frac{f(A=1\mid w)}{f(A=0\mid w)}\cd y\cd f(y,A=0,w)\,dy\,dw\\
        &= \E\bigg[\frac{\one(A=0)}{f(A=1)}\cd \frac{f(A=1\mid W)}{f(A=0\mid W)}\cd\E[Y\mid A,W] \bigg] 
    \end{split}
\end{equation}
Consequently, the Riesz representer is \[\alpha(A,W) = \frac{\one(A=0)}{f(A=1)}\cd \frac{f(A=1\mid W)}{f(A=0\mid W)}.\] And, for the NDE (Example \ref{nde}):
\begin{equation}\label{eq:ipw:nde}
    \begin{split}
        \theta(a')&= \E[\E\{\E[Y\mid A=a', M, W]\mid A=0, W]\}] \\
            &= \iiint y \cd f(y\mid A=a',m,w)\cd f(m\mid A=0,w)\cd f(w)\,dy\,dm\,dw \\
            &=\iiint y \cd \frac{f(y,A=a',m,w)}{f(m,A=a',w)}\cd \frac{f(m,A=0,w)}{f(A=0,w)}\cd f(w)\,dy\,dm\,dw \\
            &= \iiint y \cd \frac{f(y,A=a',m,w)}{f(m\mid A=a',w)f(A=a'\mid w)f(w)}\cd \frac{f(m\mid A=0,w)f(A=0\mid w)f(w)}{f(A=0\mid w)f(w)}\cd f(w)\,dy\,dm\,dw \\
            &=\iiint \frac{y\cd f(y,A=a',m,w)}{f(A=a'\mid w)}\cd \frac{f(m\mid A=0,w)}{f(m\mid A=a', w)}\,dy\,dm\,dw \\
            &= \E\bigg[\frac{\one(A=a')}{f(A=a'\mid W)}\frac{f(M\mid A=0,W)}{f(M\mid A=a',W)}\cd \E[Y\mid A,M,W]\bigg].
    \end{split}
\end{equation}
By linearity of expectation, it follows that the Riesz representer for the NDE is then the difference of the Riesz representers for $\theta(1)$ and $\theta(0)$.
\begin{equation}
    \alpha(A,M,W) = \frac{\one(A=a')}{f(A=a'\mid W)}\frac{f(M\mid A=0,W)}{f(M\mid A=a',W)} - \frac{\one(A=0)}{f(A=0\mid W)}.
\end{equation}

\section{Using Riesz representers to derive the EIF}\label{section:eif}

Parametric models are popular because we can use classical statistical theory to quantify uncertainty around estimates of the statistical estimand (e.g., calculate confidence intervals). However, the estimates produced by parametric models may be biased if those models are misspecified, which in practice, is likely. Instead, semiparametric estimators can leverage flexible machine learning algorithms in model fitting. Such commonly used estimators include AIPW,\citep{robins_estimation_1994} targeted minimum-loss based estimation (TMLE),\citep{van2006targeted} and double/debiased machine learning (DML).\citep{chernozhukov2018double} As stated in \S\ref{section:intro}, these estimators are based on the EIF of the statistical estimand. A full review of the EIF is outside the scope of this work, as it has been extensively covered elsewhere.\citep{hampel1974influence,fisher2021visually,Hines_Dukes_Diaz-Ordaz_Vansteelandt_2022,renson2025pulling} Briefly, the EIF is a mean-zero derivative that quantifies how sensitive the statistical estimand is to changes in the underlying data distribution. The variance of the EIF also represents a lower bound on the variance of an unbiased estimator of the statistical estimand. Consequently, this means that estimands based on EIFs are efficient in the sense that they asymptotically achieve this variance bound. 

The Riesz representer of a statistical estimand has a convenient relationship to the EIF in that the existence of the Riesz representer means that the EIF of the statistical estimand is finite\citep{chernozhukov2022automatic} (i.e., the estimand is referred to as being pathwise differentiable). In addition, if the statistical estimand can be expressed as sequential mappings of iterated expectations (such as the example estimands we consider in this paper), then one can use the Riesz representation of the statistical estimand as part of a simple recursive algorithm for deriving the semiparametric EIF.\citep{newey1994asymptotic,chernozhukov2022automatic} The algorithm is as follows: 

\begin{alg}\label{eifalg}
    Let $\theta$ be the parameter of interest that can be defined in terms of sequential expectations, $K$, which is equal to the number of expectations, and $\m_{K+1}(\Q_{K+1};X) = Y$. Recursively, for $k = K, ..., 1$:
    \begin{enumerate}
        \item Let $\Q_k(x)$ represent the innermost expectation and $\m_k(\Q_k;X)$ be the mapping of that expectation.
        \item Derive the Riesz representer $\alpha_k(X)$ of the parameter such that
        \[\theta = \cd\E[\alpha_k(X)\cd\m_{k+1}(\Q_{k+1};X)] = \cd\E[\alpha_k(X)\Q_k(X)].\] 
        \item Let
            $$
            D_k = \begin{cases}
                \alpha_k(X)\{\m_{k+1}(\Q_{k+1};X) - \Q_k(X)\} &\text{ if } k > 1,\\
                \alpha_k(X)\{\m_{k+1}(\Q_{k+1};X) - \theta\} &\text{ if } k = 1\\
            \end{cases}
            $$
        \item If $k > 1$, treat $\m_k(\Q_k;X)$ as a random variable and recursively update the estimand definition. Otherwise, exit the recursion.
    \end{enumerate}
    The nonparametric EIF of $\theta$ is
    \[\phi(O; \eta) = \cd\sum_{k=1}^K D_k\]
    where $\eta = (\alpha_1(X), \Q_k(X), ..., \alpha_K(X), \Q_K(X))$ denotes the nuisance parameters.
\end{alg}

We can apply the algorithm to derive the EIF of the example estimands. We omit applying the algorithm to Example \ref{rct} for brevity.  

\noindent \textbf{Example 2} \textit{ATE} As a reminder, the ATE is given by Eq. \ref{eq:ate}, which has $K = 2$ regressions.
\begingroup 
\renewcommand{\theenumi}{$k=$ \arabic{enumi}}
\begin{etaremune}[itemsep = 0em,leftmargin = 0.5in]
    \item \,
        \begin{enumerate}[label=Step. \arabic*., itemsep = -0.5em]
            \item \textit{Define the regression and mapping.}
            	$$\Q_2(a, w) = \E[Y\mid A=a,W=w]\text{ and }\m_2(\Q_2;W) = \Q_2(1, W) - \Q_2(0, W)$$
            \item \textit{Derive the Riesz representer.} From Eq. \ref{eq:ipw:ate}, we know $$\alpha_2(A,W) = \frac{\one(A=1)}{f(A=1\mid W)} - \frac{\one(A=0)}{f(A=0\mid W)}$$
            \item \textit{Derive the EIF component.} We have, $$D_2 = \bigg\{\frac{\one(A=1)}{f(A=1\mid W)} - \frac{\one(A=0)}{f(A=0\mid W)} \bigg\}\{Y - \E[Y\mid A,W]\}$$
            \item \textit{Update the estimand definition.} Recursively updating the estimand yields $$\theta = \E[\m_2(\Q_2;W)].$$
        \end{enumerate}
    \item \,
    	\begin{enumerate}[label=Step. \arabic*., itemsep = -0.5em]
	 	\item $\Q_1(a) = \E[\m_2(\Q_2;W)]$ and $\m_1(\Q_1) = \Q_1$.
		\item In this case, the Riesz representer is just $1$. So, $\alpha_1 = 1$. 
		\item $$D_1 = \E[Y\mid A=1,W] -  \E[Y\mid A=0,W]  - \theta.$$
		\item Because $k =1$, we exit the recursion.
	\end{enumerate}
\end{etaremune}
\endgroup 
\noindent Summing the $D_k$ terms, the nonparametric EIF of the ATE\citep{vanDerLaan2011TargetedLearning} is given by
\begin{equation}
    \begin{split}
        \phi(O;\eta) = &\bigg\{\frac{\one(A=1)}{f(A=1\mid W)} - \frac{\one(A=0)}{f(A=0\mid W)}\bigg\}\{Y - \E[Y\mid A,W]\}\\ &+ \E[Y\mid A=1,W] - \E[Y\mid A=0,W] - \theta.
    \end{split}
\end{equation}

\noindent \textbf{Example 3} \textit{ATT} The mean outcome under control for the treated is identified as Eq. \ref{eq:att}. Again, the number of regressions is $K = 2$.
\begingroup 
\renewcommand{\theenumi}{$k=$ \arabic{enumi}}
\begin{etaremune}[itemsep = 0em,leftmargin = 0.5in]
    \item \,
        \begin{enumerate}[label=Step. \arabic*., itemsep = -0.5em]
            \item \textit{Define the regression and mapping.} $$\Q_2(a,w) = \E[Y\mid A=a,W=w]\text{ and }\m_2(\Q_2;W) = \Q_2(0,W).$$
            \item \textit{Derive the Riesz representer.} From Eq. \ref{eq:ipw:att}, we know $$\alpha_2(A,W) = \frac{\one(A=0)}{f(A=1)}\frac{f(A=1\mid W)}{f(A=0\mid W)}.$$
            \item \textit{Derive the EIF component.} We have, $$D_2 = \frac{\one(A=0)}{f(A=1)}\frac{f(A=1\mid W)}{f(A=0\mid W)}\{Y - \E[Y\mid A,W]\}.$$
            \item \textit{Update the estimand definition.} Performing the recursion yields $$\theta = \E[\m_2(\Q_2;W)\mid A= 1].$$
        \end{enumerate}
    \item \,
    	\begin{enumerate}[label=Step. \arabic*., itemsep = -0.5em]
	 	\item The regression and mapping is now $$\Q_1(a) = \E[\m_2(\Q_2;W)\mid A=a]\text{ and }\m_1(\Q_1) = \Q_1(1).$$
		\item The current estimand definition is similar to that of Example \ref{rct}. Therefore, from Eq. \ref{eq:ipw:rct}, we have $$\alpha_1(A) = \frac{\one(A = 1)}{f(A = 1)}.$$
		\item $$D_1 = \frac{\one(A = 1)}{f(A = 1)}\{\E[Y\mid A=0,W] - \theta\}.$$
		\item Exit the recursion.
	\end{enumerate}
\end{etaremune}
\endgroup 
\noindent The nonparametric EIF of the mean outcome under control for the treated is then
\begin{equation}
    \begin{split}
        \phi(O;\eta) &= \frac{\one(A=0)}{f(A=1)}\frac{f(A=1\mid W)}{f(A=0\mid W)}\{Y - E[Y\mid A, W]\} \\&+ \frac{\one(A = 1)}{f(A=1)}\big\{\E[Y\mid A =0, W] - \theta\big\}.\citep{vanDerLaan2011TargetedLearning}
    \end{split}
\end{equation}

\noindent \textbf{Example 4} \textit{NDE} Recall that the NDE is a contrast of Eq. \ref{eq:nde}. Here, $K = 3$. 
\begingroup 
\renewcommand{\theenumi}{$k=$ \arabic{enumi}}
\begin{etaremune}[itemsep = 0em,leftmargin = 0.5in]
    \item \,
        \begin{enumerate}[label=Step. \arabic*., itemsep = -0.5em]
            \item \textit{Define the regression and mapping.} $$\Q_3(a,m,w) = \E[Y\mid A=a,M=m,W=w]\text{ and }\m_3(\Q_3;M,W) = \Q_3(a', M, W).$$
            \item \textit{Derive the Riesz representer.} From Eq. \ref{eq:ipw:nde}, we have $$\alpha_3(A,M,W) = \frac{\one(A=a')}{f(A=a'\mid W)}\frac{f(M\mid A=0,W)}{f(M\mid A=a',W)}.$$
            \item \textit{Derive the EIF component.} $$D_3 =\frac{\one(A=a')}{f(A=a'\mid W)}\frac{f(M\mid A=0,W)}{f(M\mid A=a',W)}\{Y - \E[Y\mid A,M,W]\}.$$
            \item \textit{Update the estimand definition.} $\theta(a') =\E\{\E[\m_3(\Q_3;M,W)\mid A=0, W]\}.$
        \end{enumerate}
    \item \,
    	\begin{enumerate}[label=Step. \arabic*., itemsep = -0.5em]
	 	\item The regression and mapping is now $$\Q_2(a,w) = \E[\m_3(\Q_3;M,W)\mid A=a,W=w]\text{ and }\m_2(\Q_2;W) = \Q_2(0,W).$$
		\item The current estimand definition is similar to the parameter $\E\{\E[Y\mid A=0,W]\}$, which is a component of the ATE. Therefore, from Eq. \ref{eq:ipw:ate}, we have $$\alpha_2(A,W) = \frac{\one(A=0)}{f(A=0\mid W)}.$$
		\item $$D_2 =\frac{\one(A=0)}{f(A=0\mid W)}\big\{\E[Y\mid A=a',M,W] - \E\{E[Y\mid A=a',M,W]\mid A,W \}\big\}.$$
		\item $\theta(a') = \E[\m_2(\Q_2;W)]$.
	\end{enumerate}
    \item \,
    	\begin{enumerate}[label=Step. \arabic*., itemsep = -0.5em]
	 	\item $\Q_1 = \E[\m_2(\Q_2;W)]$ and $\m_1(\Q_1) = \Q_1$.
		\item $\alpha_1 = 1$.
		\item $$D_1 =\E\{\E[Y\mid A=a',M,W]\mid A=0,W\} - \theta(a')$$
		\item Exit the recursion.
	\end{enumerate}
\end{etaremune}
\endgroup 
\noindent Summing the $D_k$ terms yields
\begin{equation}\label{eq:eif:nde}
    \begin{aligned}
        \phi(O; \eta) &= \frac{\one(A=a')}{f(A=a'\mid W)}\frac{f(M\mid A=0,W)}{f(M\mid A=a',W)}\{Y - \E[Y\mid A, M, W]\} \\
        &+ \frac{\one(A=0)}{f(A=0\mid W)}\{\E[Y\mid A=a',M,W] - \E[\E[Y\mid A=a',M,W]\mid A,W]\} \\
        &+\E\{\E[Y\mid A=a',M,W]\mid A=0,W\} - \theta(a')
    \end{aligned}
\end{equation}
as the nonparametric EIF of Eq. \ref{eq:nde}.\citep{tchetgen2012semiparametric} 

\section{Riesz regression}

So far, we have used Riesz representers as a convenient theoretical tool for deriving EIFs of certain reweighting (IPW-style) estimands. However, the real strength of this framework lies in the development of so-called Riesz regression methods, which can directly learn Riesz representers from data. As we briefly discussed, semiparametric estimators that use machine learning are based on the EIF of the statistical estimand. For instance, TMLE iteratively updates initial estimates of nuisance parameters until the EIF has an expected value of zero. These procedures, therefore, require generating estimates of each of the nuisance parameters that appear in the EIF. Some nuisance parameters are easy to estimate with off-the-shelf algorithms. For example, in the case of the ATE, we can estimate the propensity score using any machine learning algorithm capable of binary classification. For many other statistical estimands, however, estimating one or more of the EIF's nuisance parameters is hard. Consider Eq. \ref{eq:nde}; its EIF (Eq. \ref{eq:eif:nde}) includes a conditional density ratio of the mediator: $f(M \mid A=0, W)\cd f(M \mid A=1, W)^{-1}$. If the mediator is high-dimensional (e.g., consisting of many variables), then this density ratio can be challenging to estimate using off-the-shelf machine learning algorithms. (We note that under this simple data generating mechanism, one could re-parameterize the EIF to include only easy-to-estimate nuisance parameters,\cite{Diaz2021-gi} but such reparameterization is not possible under more complex---and likely more realistic---data generating mechanisms that involve non-discrete (e.g., continuous or multivariate) intermediate confounders and non-discrete mediators.\cite{Rudolph2024-mx} In these more complex settings where nuisance parameters involving multivariate or continuous conditional densities cannot be re-parameterized away, Riesz regression may be a more practical solution.)

As an alternative, we can directly estimate Riesz representers with the so-called Riesz loss function.\citep{chernozhukov2021automatic} An estimator that minimizes the empirical expectation of the Riesz loss is referred to as a Riesz regression. In what follows, $f_\theta(X)$ denotes some model with parameters $\theta$ that takes as input the variables $X$ and returns a numeric vector. For example, $f_\theta$ could be a neural network with node weights $\theta$.\citep{chernozhukov2022riesznet} The Riesz representer $\alpha(X)$ can be directly estimated by minimizing the empirical mean-squared error loss within some function class $\mathcal{F}_n$: 
\begin{equation}\label{rieszloss}
    \begin{split}
        \hat{\alpha}(X) &= \argmin_{f_\theta \in \mathcal{F}_n}\E[(\rr - f_\theta(X))^2]\\
        &= \argmin_{f_\theta \in \mathcal{F}_n}\E[\rr^2 - 2\rr f_\theta(X) + f_\theta(X)^2] \\
        &= \argmin_{f_\theta \in \mathcal{F}_n}\E[\rr^2] - 2\E[\rr f_\theta(X)] + \E[f_\theta(X)^2]\\
        &= \argmin_{f_\theta \in \mathcal{F}_n}\E[f_\theta(X)^2] - 2\E[\rr f_\theta(X)]\\
        &= \argmin_{f_\theta \in \mathcal{F}_n}\E[f_\theta(X)^2 - 2\m_k(f_\theta;X)].
    \end{split}
\end{equation}
The second equality follows from squaring the binomial, the fourth by recognizing that $\E[\rr^2]$ is a constant that doesn't depend on $f_\theta(X)$, and the fifth by the Riesz representation theorem. Here, to clarify, the notation $\m_k(f_\theta; X)$ signifies that the mapping $\m_k$ is applied to the function $f_\theta$ instead of the regression $\Q_k$. For example, when $\m_2(\Q_k;W) = \Q_2(1, W) - \Q_2(0, W)$ (see Example \ref{ate}), the expression $\m_2(f_\theta;W)$ becomes $f_\theta(1, W) - f_\theta(0, W)$. This effectively calculates the difference in $f_\theta$ when A is set to 1 versus when A is set to 0. Similarly, for Example \ref{att}, $\m_2(f_\theta;W) = f_\theta(0, W)$; and $\m_2(f_\theta;M,W) = f_\theta(a', M, W)$ from Example \ref{nde}.

For instance, a semiparametric estimator for Eq. \ref{eq:nde} based on the EIF (Eq. \ref{eq:eif:nde}) requires estimating two Riesz representers. While we previously derived their analytical form, let's proceed as if they were unknown to us. Beginning with $\alpha_2(A,W)$, we have
\begin{equation}\label{eq:rieszloss:alpha2}
\begin{aligned}
    \hat{\alpha}_2(A,W) 
        &= \argmin_{f_\theta \in \mathcal{F}_n} \E[f_\theta(A,W)^2 - 2\alpha_2(A,W)f_\theta(A,W)]\\
        &= \argmin_{f_\theta \in \mathcal{F}_n} \E[f_\theta(A,W)^2 - 2\m_2(f_\theta;W)]\\
        &= \argmin_{f_\theta \in \mathcal{F}_n} \E[f_\theta(A,W)^2 - 2f_\theta(0; W)].
\end{aligned}        
\end{equation}
Similarly,
\[\hat{\alpha}_3(A,M,W) = \argmin_{f_\theta \in \mathcal{F}_n} \E[f_\theta(A,M,W)^2 - 2\alpha_3(A,M,W)f_\theta(A,M,W)].\]
From Theorem \ref{riesz} and Algorithm \ref{eifalg}, we have the following set of equalities
\begin{equation}\label{eq:equalities:nde}
\begin{split}
  \theta(a')
  &= \E\Big[\E\big\{\E[Y \mid A=a', M, W] \mid A=0, W\big\}\Big] \\
  &= \E\big[\alpha_3(A,M,W) \cd Y\big] \\
  &= \E\big[\alpha_3(A,M,W)\,\Q_3(A,M,W)\big] \\
  &= \E\big\{\E\big[\m_3(\Q_3;M,W) \mid A=0, W\big]\big\} \\
  &= \E\big[\alpha_2(A,W)\,\m_3(\Q_3;M,W)\big] \\
  &= \E\big[\alpha_2(A,W)\,\Q_2(A,W)\big] \\
  &= \E\big[\m_2(\Q_2;W)\big].
\end{split}
\end{equation}
This implies that $\E[\alpha_3(A,M,W)f_\theta(A,M,W)] = \E[\alpha_2(A,W)\m_3(f_\theta;M,W)]$. Performing this substitution yields
\begin{equation}\label{eq:rieszloss:alpha3}
\begin{aligned}
    \hat{\alpha}_3(A,M,W) 
        &= \argmin_{f_\theta \in \mathcal{F}_n} \E[f_\theta(A,M,W)^2 - 2\hat{\alpha}_2(A,W)\m_3(f_\theta;M,W)] \\
        &= \argmin_{f_\theta \in \mathcal{F}_n} \E[f_\theta(A,M,W)^2 - 2\hat{\alpha}_2(A,W)f_\theta(a',M,W)],
\end{aligned}        
\end{equation}
where the last equality follows from the fact that $\m_3(\Q_3;M,W) = \Q_3(a',M,W)$. This indicates that we can directly estimate the Riesz representers for Eq. \ref{eq:nde} sequentially by first estimating $\alpha_2$ and then using this estimate as weights in the Riesz loss for learning $\alpha_3$. 

In practice, Riesz regression has been performed using deep learning,\citep{chernozhukov2022riesznet} random forests,\citep{chernozhukov2022riesznet} and gradient boosting.\citep{lee2025rieszboost} We provide example code of the deep learning approach using \texttt{R} in the eAppendix. The key insight to recognize is that these loss functions don't require knowing the closed form expression of the Riesz representers. For this reason, estimators that rely on this strategy are often referred to as automatic.\citep{chernozhukov2022automatic} 

\section{Discussion}

Our goal was to provide an accessible introduction to the Riesz representation theorem in the context of semiparametric estimators based on the efficient influence function. Riesz regression can meaningfully improve our ability to construct semiparametric estimators that depend on nuisance parameters that would otherwise be extremely difficult to estimate. For example, \citeauthor{liu2024general} use sequential Riesz regressions via deep learning for a one-step estimator of mediation effects with high-dimensional mediators and mediator-outcome confounders. Similarly, \citeauthor{susmann2024longitudinal} use Riesz regression in a TMLE for estimating a generalization of the ATT in longitudinal data. 

Although the adoption of Riesz regression is accelerating, practical barriers to implementation remain. Currently, there are limited off-the-shelf software implementations that have been thoroughly tested for minimizing the Riesz loss. Consequently, using Riesz regression may require writing bespoke code and, thus, a greater understanding of tools like deep learning or boosting. Additionally, simulation evidence for the benefit of using Riesz regression has, so far, primarily come from papers that propose new methods, and more independent assessments of its performance are required. For these reasons, we recommend using Riesz regression carefully and comparing results with more established methods where possible.
Although we considered example statistical parameters 1-4 in this paper, it would only make practical sense to use Riesz regression in Example 4 in cases with multivariate or continuous mediators. However, one could still avoid Riesz regression in Example 4 through a re-parameterization of the mediator conditional densities.\citep{Diaz2021-gi,Rudolph2024-mx} In Examples 1-3, each Riesz representer would be easily estimated with standard off-the-shelf software. 

We sought to address other practical barriers to implementing Riesz regression: lack of understanding due to the current literature being too technical for applied audiences, and lack of clear motivation for how Riesz regression could address real-world challenges. With those goals in mind, we end with a one-sentence summary: when you encounter the term ``Riesz representer'' in the context of semiparametric statistical estimators, simply think ``weights in a reweighting estimator.''

\newpage

\bibliography{lib}

\appendix

\definecolor{codegreen}{rgb}{0,0.6,0}
\definecolor{codegray}{rgb}{0.5,0.5,0.5}
\definecolor{codepurple}{rgb}{0.58,0,0.82}
\definecolor{backcolour}{rgb}{0.95,0.95,0.92}

\lstdefinestyle{mystyle}{
    backgroundcolor=\color{backcolour},   
    commentstyle=\color{codegreen},
    keywordstyle=\color{magenta},
    numberstyle=\tiny\color{codegray},
    stringstyle=\color{codepurple},
    basicstyle=\ttfamily,
    showstringspaces=false,
    numbers=left
}

\lstset{style=mystyle}

\section*{eAppendix}

\begin{definition}[Linear mapping]
An mapping, $\m_k$, is linear if $\m_k(c_1\cdot\Q_{k,1} + c_2\cdot\Q_{k,2}; X) = c_1\cdot\m_k(\Q_{k,1};X) + c_2\cdot\m_k(\Q_{k,2};X)$ for any functions $\Q_{k,1}$, $\Q_{k,2}$ and scalars $c_1$, $c_2$.
\end{definition}

\begin{definition}[Bounded linear functional]
A linear functional $\m_1(\Q_1)$ is bounded if there is a constant (i.e., a scalar) $C > 0$ where $|\m_1(\Q_1)| \leq C\sqrt{\E[\Q_1(X)^2]}$ for all functions $\Q_1(X)$ where $\E[\Q_1(X)^2]<\infty$. In words, for functions $\Q_1(X)$ in the the $L^2_P$ space, a functional is bounded if the magnitude (i.e. absolute value) its value is no greater than the $L^2_P$ norm of $\Q_1(X)$, the input, up to a constant.
\end{definition}

\subsection*{Example of Riesz regression in \texttt{R}}

As a simple example, we consider directly estimating the Riesz representer $\alpha_3(A,M,W)$ for the parameter
\[\theta = \E[\E\{\E[Y\mid A=1, M, W]\mid A=0, W]\}].\]
As was shown in the main text, we can directly estimate $\alpha_3$ with Riesz regression by minimizing the Riesz loss:
$$
\alpha_3(A,M,W) = \argmin_{f_\theta} \E[f_\theta(A,M,W)^2 - 2\alpha_2(A,W)f_\theta(a',M,W)],
$$
where $\alpha_2(A,W) = \frac{\one(A=1)}{f(A=1\mid W)}$. We use simulated data described by the following data generating mechanism:
\begin{align*}
    P(W = 1) &= 0.4 \\
    P(A = 1) &= 0.3 + 0.3W \\
    M &\sim \mathcal{N}(0.6 + 0.05A -0.3W, 1)\\
    P(Y = 1\mid A,M,W) &= \text{expit}(-\text{log}(5) + \text{log}(2)A + \text{log}(3)M - \text{log}(1.2)W).
\end{align*}
As our Riesz regression, we use a multi-layer perceptron (MLP) neural network architecture.\citep{goodfellow2016deep} This neural network consists of $K=2$ hidden linear layers, each with 4 neurons and a Rectified Linear Unit (ReLU) activation function.\citep{nair2010rectified} The model is trained using gradient descent with the Adam optimizer.\citep{kingma2014adam} We implement the neural network in \texttt{R}\citep{Rlang} with the \textit{torch} package.\citep{torchR} We emphasize that this is a toy example. Importantly, the following approach should be cross-fit to avoid invoking a Donsker class assumption.\citep{kennedy2024semiparametric} Other optimizations could include: early stopping based on loss in a validation set, batch training, and dynamic learning rate adjustment.

\vspace{1em}
\begin{lstlisting}[language=R]
library(torch)
library(progressr)

set.seed(56)
n <- 5000
W <- rbinom(n, 1, 0.4)
A <- rbinom(n, 1, 0.3 + 0.3*W)
M <- 0.6 + 0.05*A - 0.3*W + rnorm(n)
Y <- rbinom(n, 1, plogis(-log(5) + log(2)*A + 
                           log(3)*M - log(1.2)*W))
example_data <- data.frame(W, A, M, Y)

aprime <- 1

# Estimate \alpha_2(A,W) using a GLM
propensity_fit <- glm(A ~ W, data = example_data, 
                      family = "binomial")
propensity <- predict(propensity_fit, type = "response")
alpha_2 <- torch_tensor(as.numeric(A == aprime) / propensity)

# We need to convert the data to torch tensors.
# For 'mapped', we apply the function d.
observed <- torch_tensor(
  as.matrix(example_data[, c("A", "M", "W")])
)
mapped <- observed$clone()
mapped[, 1] <- aprime

neurons <- 4

# We'll use a multi-layer perceptron with two hidden layers 
# (each with 4 nuerons) and the ReLU activation function. 
# The network takes as input (A,W) and outputs a single value.
f_theta <- nn_sequential(
  nn_linear(ncol(observed), neurons),
  nn_relu(),
  nn_linear(neurons, neurons),
  nn_relu(),
  nn_linear(neurons, neurons),
  nn_relu(),
  nn_linear(neurons, 1),
  nn_relu()
)

# Set the learning rate to 0.01
optimizer <- optim_adam(params = f_theta$parameters, lr = 0.01)

# Train the network over 200 epochs.
epochs <- 200

torch_manual_seed(634523)

with_progress({
  # Create a progress bar to visualize training progress
  p <- progressor(steps = epochs)
  
  for (epoch in 1:epochs) {
    # Riesz loss
    riez_loss <- (f_theta(observed)$pow(2) - 
                    (2 * alpha_2 * f_theta(mapped)))$mean()
    
    if (epoch %% 10 == 0)
      cat("Epoch: ", epoch, " Loss: ", riez_loss$item(), "\n")
    
    optimizer$zero_grad()
    riez_loss$backward()
    
    optimizer$step()
    p()
  }
})
#> Epoch:  10  Loss:  -0.4690982                                                   
#> Epoch:  20  Loss:  -0.8258173                                                   
#> Epoch:  30  Loss:  -1.060881                                                    
#> Epoch:  40  Loss:  -1.134421                                                    
#> Epoch:  50  Loss:  -1.401238                                                    
#> Epoch:  60  Loss:  -2.075226                                                    
#> Epoch:  70  Loss:  -2.614951                                                    
#> Epoch:  80  Loss:  -2.760663                                                    
#> Epoch:  90  Loss:  -2.761145                                                    
#> Epoch:  100  Loss:  -2.766564                                                   
#> Epoch:  110  Loss:  -2.769224                                                   
#> Epoch:  120  Loss:  -2.769007                                                   
#> Epoch:  130  Loss:  -2.769578                                                   
#> Epoch:  140  Loss:  -2.769637                                                   
#> Epoch:  150  Loss:  -2.769673                                                   
#> Epoch:  160  Loss:  -2.769707                                                   
#> Epoch:  170  Loss:  -2.76972                                                    
#> Epoch:  180  Loss:  -2.76973                                                    
#> Epoch:  190  Loss:  -2.769735                                                   
#> Epoch:  200  Loss:  -2.769739

f_theta$eval()

# Return the estimated Riesz representers
alpha_3 <- as.numeric(f_theta(observed))

# Compare to the truth
true_alpha_2 <- as.numeric(A == aprime) / (0.3 + 0.3*W)
true_alpha_3 <- true_alpha_2 * 
  (dnorm(M, 0.6 - 0.3*W) / dnorm(0.6 + 0.05 - 0.3*W))

# MSE
mean((alpha_3 - true_alpha_3)^2)
# [1] 0.3515668
\end{lstlisting}

\end{document}